\documentclass[12pt]{amsart} 
\usepackage{amssymb}
\usepackage{latexsym}
\usepackage{amsmath}
\usepackage{amsthm}
\usepackage{amsfonts}

\usepackage{soul}
\usepackage[usenames,dvipsnames]{xcolor}

\usepackage{doi}

\usepackage{mathabx}

\newtheorem*{nonumbertheorem}{Theorem}

\newtheorem*{nonumbercorollary}{Corollary}

\theoremstyle{definition}

\newtheorem*{acknowledgement}{Acknowledgement}
\theoremstyle{remark}

\usepackage{color}

\newcommand{\R}{\mathbb R}
\newcommand{\N}{\mathbb N}
\newcommand{\C}{\mathbb C}

\title
[Non-extendability  of the finite Hilbert transform]
{Non-extendability  of the finite \\  Hilbert transform}

\author[G.P.  Curbera]{Guillermo P. Curbera}
\address{Facultad de Matem\'aticas \& IMUS,
Universidad de Sevilla, 
Calle Tarfia s/n,  Sevilla 41012, Spain}
\email{curbera@us.es}

\author[S. Okada]{Susumu Okada} 
\address{School of Natural Sciences (Maths/ Physics), University of Tasmania, 
Private Bag 37, Hobart, Tas. 7001, Australia}
\email{susumu.okada@utas.edu.au}

\author[W.J. Ricker]{Werner J. Ricker}
\address{Math.--Geogr.\  Fakult\"at, Katholische Universit\"at
Eichst\"att--Ingolstadt, D--85072 Eichst\"att, Germany}
\email{werner.ricker@ku.de}

\thanks{The first author acknowledges the support  of 
of PGC2018-096504-B-C31, FQM-262 and Feder-US-1254600 (Spain).}

\date{\today}



\begin{document}

\begin{abstract}
The finite Hilbert transform $T\colon X\to X$  acts continuously 
on every rearrangement invariant space $X$ on  $(-1,1)$ having non-trivial Boyd indices.
It is proved that $T$ cannot be 
further extended, whilst still taking its values in $X$, to any larger domain space.
That is, $T\colon X\to X$ is  already optimally defined. 
\end{abstract}

\maketitle


\section{Introduction and main result}
\label{S1}


The finite Hilbert transform $T(f)$ 
of $f\in L^1(-1,1)$ is the well known principal value integral
\begin{equation*}
(T(f))(t)=\lim_{\varepsilon\to0^+} \frac{1}{\pi}
\left(\int_{-1}^{t-\varepsilon}+\int_{t+\varepsilon}^1\right) \frac{f(x)}{x-t}\,dx ,
\end{equation*}
which exists for a.e.\ $t\in(-1,1)$ and is a measurable function. 
It  has important applications to aerodynamics  and elasticity
via the  airfoil  equation, 
\cite{cheng-rott}, \cite{muskhelishvili}, \cite{reissner},  \cite{tricomi-1},
\cite{tricomi}, and
to problems arising in image reconstruction; 
see, for example, \cite{katsevich-tovbis}, \cite{sidky-etal}.
We refer to \cite{duduchava}, \cite{gakhov}, \cite{gohberg-krupnik-1},
\cite{gohberg-krupnik-2},  \cite{mikhlin-prossdorf}
where one-dimensional singular integral operators 
closely related to the finite Hilbert transform are studied in great detail

For each $1<p<\infty$
the classical linear operator $f\mapsto T(f)$ maps $L^p(-1,1)$ continuously  into 
itself; denote this operator by $T_p$.  Tricomi showed that 
$T_p$ is a Fredholm operator and  exhibited  
inversion formulae, \cite{tricomi-1}, except for the case when $p=2$, 
\cite[\S4.3]{tricomi}  (see also \cite[Ch. 11]{king}, \cite[Ch. 14.4-3]{polyanin-manzhirov} and the references therein). For $T_2$ the situation is significantly 
different, as already pointed 
out somewhat earlier in \cite[p.44]{sohngen}. Partial operator theoretic results 
for $T_2$ on $L^2(-1,1)$ 
were obtained by Okada and  Elliott, \cite{okada-elliot}; see also the references.

In \cite{curbera-okada-ricker} the finite 
Hilbert transform $T$  was studied when acting on 
suitable rearrangement invariant (r.i., in short)  spaces 
$X$ on $(-1,1)$; see below for  the relevant definitions. 
Actually, $T$ acts continuously on $X$
(denote this operator by $T_X$) precisely when the Boyd indices of $X$ 
are non-trivial, that is, when
$0<\underline{\alpha}_X\le\overline{\alpha}_X<1$; 
see \cite[pp.170--171]{krein-petunin-semenov}. This class of 
r.i.\ spaces is the largest and most  adequate replacement  for the 
$L^p$-spaces when undertaking a further
study of the finite Hilbert transform $T$. 
This is due to  two critical facts: that 
$T\colon X\to X$ is injective if and only if the function $1/\sqrt{1-x^2}\notin X$,  
and  that $T\colon X\to X$ has non-dense range  
if and only if $1/\sqrt{1-x^2}$ belongs to the associate space $X'$ of $X$
 (whenever  $X$ is separable). 
In terms of r.i.\ spaces the previous conditions can be phrased as follows:
$T\colon X\to X$ \textit{is injective if and only if $L^{2,\infty}(-1,1)\not\subseteq X$} 
and \textit{$T\colon X\to X$ has a 
non-dense range if and only if $X \subseteq L^{2,1}(-1,1)$}  (for $X$ separable). 
Here $L^{2,1}(-1,1)$ and $L^{2,\infty}(-1,1)$ are 
the usual Lorentz spaces.

Various types of inversion  results of Tricomi  for 
the operator $T_p$ (when $1<p<2$ and $2<p<\infty$) have been
extended to $T_X$ whenever the Boyd indices of $X$  satisfy
the condition $0<\underline{\alpha}_X\le\overline{\alpha}_X<1/2$ or
$1/2<\underline{\alpha}_X\le\overline{\alpha}_X<1$; 
see  \cite[Theorems 3.2 and 3.3]{curbera-okada-ricker}.
Moreover, $T$ is necessarily a Fredholm operator
in such r.i.\ spaces, \cite[Remark 3.4]{curbera-okada-ricker}.
Results  of this kind admit the possibility for a refinement of the solution of the 
airfoil equation; see \cite[Corollary 3.5]{curbera-okada-ricker}. 
Additional operator theoretic results concerning $T_X$ in r.i.\ 
spaces $X$ occur in the recent article \cite{curbera-okada-ricker2} 
(e.g., compactness, order boundedness, integral representation,
etc.).

An important problem is the possibility of extending the domain
of $T_p$, with $T_p$  still maintaining its values in $L^p(-1,1)$. 
It was shown in \cite[Example 4.21]{okada-ricker-sanchez}, 
for all  $1<p<\infty$ with $p\not=2$, that   there is \textit{no} 
larger Banach function space (B.f.s.\ in short)   containing
$L^p(-1,1)$ such that $T_p$ has an $L^p(-1,1)$-valued continuous extension to this space.
This result was generalized  in \cite[Theorem 4.7]{curbera-okada-ricker}. 
Namely, it is not possible to  extend the finite 
Hilbert transform $T_X\colon X\to X$ for any r.i.\ space $X$ satisfying 
\begin{equation}\label{a}
0<\underline{\alpha}_X\le\overline{\alpha}_X<1/2
\quad\textrm{or}\quad
1/2<\underline{\alpha}_X\le\overline{\alpha}_X<1.
\end{equation}

The arguments used in  \cite{curbera-okada-ricker} for 
establishing the above result do not apply 
to $T_X$ for r.i.\ spaces $X$ which fail to satisfy \eqref{a}. 
In particular, they do not apply 
to $T_2\colon L^2(-1,1)\to L^2(-1,1)$. However, in \cite{curbera-okada-ricker} it was also established,
via a completely different approach, that at least $T_2$
does not have a continuous $L^2(-1,1)$-valued extension to any larger B.f.s.,
\cite[Theorem 5.3]{curbera-okada-ricker}.

Thus, the question of extendability of $T_X$ remains unanswered 
for a large sub-family of r.i.\ spaces which have non-trivial Boyd indices.
Indeed, with the exception of $X=L^2(-1,1)$, this is the case for all those r.i.\ spaces $X$
satisfying $0<\underline{\alpha}_X\le1/2\le \overline{\alpha}_X<1$. In particular, this includes all  the 
Lorentz spaces $L^{2,q}$ for $1\le q\le\infty$ with $q\not=2$.
The   proof given in \cite{curbera-okada-ricker} for  $T_2$, based on the Hilbert space structure of
$L^2(-1,1)$, is not applicable to other r.i.\ spaces of the kind just mentioned.

The aim of this note is to answer the above question for \textit{all} r.i.\ spaces $X$
on which $T_X$ is continuous, via a new and unified proof.

\begin{nonumbertheorem}
Let $X$ be a r.i.\ space on $(-1,1)$ with non-trivial Boyd indices. 
The finite Hilbert transform
$T_X\colon X\to X$ has no continuous, $X$-valued extension
to any genuinely larger B.f.s.\ containing $X$.
\end{nonumbertheorem}


\section{Preliminaries}
\label{S2}


In this paper the relevant measure space  is $(-1,1)$ 
equipped with its Borel $\sigma$-algebra $\mathcal{B}$ and  Lebesgue measure 
$m$ (restricted to $\mathcal{B}$). We  denote 
by $L^0(-1,1)=L^0$ the space (of equivalence classes) of all $\mathbb{C}$-valued
measurable functions, endowed with the topology of convergence in measure.
The space $L^p(-1,1)$ is denoted simply by $L^p$, for $1\le p\le\infty$.

A \textit{Banach function space} (B.f.s.) $X$ on  $(-1,1)$ is a
Banach space  $X\subseteq L^0$ satisfying
the ideal property, that is, $g\in X$ and $\|g\|_X\le\|f\|_X$
whenever $f\in X$, $g\in L^0$ and $|g|\le|f|$ a.e.  
The \textit{associate space} $X'$  of $X$ consists  of all
$g\in L^0$ satisfying $\int_{-1}^1|fg|<\infty$, for every
$f\in X$, equipped with the norm
$\|g\|_{X'}:=\sup\{|\int_{-1}^1fg|:\|f\|_X\le1\}$. 
The space $X'$ is a closed subspace of the Banach space dual $X^*$ of $X$. 
The space $X$ satisfies the Fatou property  if, whenever  $\{f_n\}_{n=1}^\infty\subseteq X$ satisfies
$0\le f_n\le f_{n+1}\uparrow f$ a.e.\ with $\sup_n\|f_n\|_X<\infty$,
then $f\in X$ and $\|f_n\|_X\to\|f\|_X$.   In this paper \textit{all} B.f.s.' $X$ 
are on $(-1,1)$ relative to Lebesgue measure and, as in 
\cite{bennett-sharpley},  satisfy the Fatou property.

A \textit{rearrangement invariant} (r.i.) space $X$ on $(-1,1)$ is a B.f.s.\  
such that if $g^*\le f^*$ with $f\in X$,  
then $g\in X$ and $\|g\|_X\le\|f\|_X$.
Here $f^*\colon[0,2]\to[0,\infty]$ is 
the decreasing rearrangement of $f$, that is, the
right continuous inverse of its distribution function:
$\lambda\mapsto m(\{t\in (-1,1):\,|f(t)|>\lambda\})$.
The associate space $X'$ of a r.i.\ space $X$ is again a r.i.\ space.
Every r.i.\ space $X$ on $(-1,1)$ satisfies 
$L^\infty\subseteq X\subseteq L^1$. 
Moreover, if $f\in X$ and $g\in X'$, then $fg\in L^1$ and
$\|fg\|_{L^1}\le \|f\|_X \|g\|_{X'}$, i.e., H\"older's inequality
is available.

The family of r.i.\ spaces includes many classical spaces 
appearing in analysis, in particular  the Lorentz $L^{p,q}$ spaces, 
\cite[Definition IV.4.1]{bennett-sharpley}.

Given a r.i.\ space $X$ on $(-1,1)$, due to the Luxemburg representation theorem 
there exists a r.i.\  space  $\widetilde X$ on $(0,2)$
such that $\|f\|_X=\|f^*\|_{\widetilde X}$ for $f\in X$, \cite[Theorem II.4.10]{bennett-sharpley}.
The dilation operator $E_t$ for $t>0$ is defined, for 
each $f\in \widetilde X$, by $E_t(f)(s):=f(st)$ for $0\le s\le \min\{2,1/t\}$ and zero 
for  $\min\{2,1/t\}< s\le 2$. The operator $E_t\colon \widetilde X\to \widetilde X$  is bounded 
with $\|E_{1/t}\|_{\widetilde X\to \widetilde X}\le \max\{t,1\}$. The \textit{lower} and \textit{upper 
Boyd indices} of  $X$ are defined, respectively, by
\begin{equation*}
\underline{\alpha}_X\,:=\,\sup_{0<t<1}\frac{\log \|E_{1/t}\|_{\widetilde X\to \widetilde X}}{\log t}
\;\;\mbox{and}\;\;
\overline{\alpha}_X\,:=\,\inf_{1<t<\infty}\frac{\log \|E_{1/t}\|_{\widetilde X\to \widetilde X}}{\log t} ;
\end{equation*}
see \cite{boyd} and also \cite[Definition III.5.12]{bennett-sharpley}. 
They satisfy $0\le\underline{\alpha}_X\le \overline{\alpha}_X\le1$.
Note that $\underline{\alpha}_{L^p}= \overline{\alpha}_{L^p}=1/p$.

For all of the above and further facts on r.i.\  spaces see \cite{bennett-sharpley}, 
for example.


\section{Proof of the Theorem}
\label{S3}


Given $X$, a r.i.\ space on $(-1,1)$ with non-trivial Boyd indices, consider the space
\begin{equation*}\label{tx}
[T,X]:=\big\{f\in L^1: T(h)\in X, \;\forall |h|\le|f|\big\},
\end{equation*}
which is a B.f.s.\ for the norm
\begin{equation*}\label{TX-norm}
\|f\|_{[T,X]}:=\sup_{|h|\le|f|} \|T(h)\|_X,\quad f\in[T,X].
\end{equation*}
The proof of this fact uses, in an essential way, a deep result of Talagrand concerning $L^0$-valued measures,
\cite[Proposition 4.5]{curbera-okada-ricker}. The space $[T,X]$ is the largest B.f.s.\ containing $X$   
to which $T_X\colon X\to X$ has a continuous, linear, $X$-valued extension, 
\cite[Theorem 4.6]{curbera-okada-ricker}.  In particular, $X\subseteq [T,X]$. Thus, in order to show that no genuine extension of $T_X$ is possible it suffices to show that $[T,X]\subseteq X$; see
 Theorems 4.7 and 5.3 in \cite{curbera-okada-ricker}.

Fix $N\in\N$. Given $a_1,\dots,a_N\in\C$   and  disjoint sets 
$A_1\dots,A_N$ in $\mathcal{B}$, define the simple function
$$
\phi:=\sum_{n=1}^N a_n\chi_{A_n}.
$$

On  $\Lambda:=\{1,-1\}^N$ consider the probability measure $d\sigma$,
which is the product measure  of $N$ copies of the 
uniform probability on $\{1,-1\}$. Define the bounded measurable function
$F$ on $\Lambda$ by
$$
\sigma=(\sigma_1,\dots,\sigma_N)\in\Lambda \mapsto F(\sigma)
:=\bigg\| T_X\bigg(\sum_{n=1}^N\sigma_n a_n\chi_{A_n}\bigg)\bigg\|_X.
$$

Observe, since the sets $A_1,\dots,A_N$ are pairwise disjoint,
that for every $\sigma=(\sigma_1,\dots,\sigma_N)\in\Lambda$ one has
$$
\left|\sum_{n=1}^N\sigma_na_n\chi_{A_n}\right|
=\sum_{n=1}^N|a_n|\chi_{A_n}=|\phi|,
$$
whence
\begin{equation}\label{1}
\|F\|_{L^\infty(\Lambda)} 
=\sup_{\sigma\in\Lambda}\left\|T\bigg(\sum_{n=1}^N\sigma_na_n\chi_{A_n}\bigg)\right\|_X
\le\sup_{|h|\le|\phi|} \|T(h)\|_X= \|\phi\|_{[T,X]}.
\end{equation}

On the other hand, an application of  Fubini's theorem yields
\begin{align}\label{2}
\|F\|_{L^\infty(\Lambda)}
 & \ge \|F\|_{L^1(\Lambda)}\nonumber
\\ &=   \int_{\Lambda}|F(\sigma)|\,d\sigma \nonumber
\\ &=  \int_{\Lambda}\bigg\| \sum_{n=1}^N\sigma_n a_nT\left(\chi_{A_n}\right)\bigg\|_X\,d\sigma \nonumber
\\ &=  \int_{\Lambda} \bigg(\sup_{\|g\|_{X'}=1} \int_{-1}^1 |g(t)| 
 \bigg| \sum_{n=1}^N\sigma_n a_nT\left(\chi_{A_n}\right)(t)\bigg|\,dt\bigg)d\sigma
 \nonumber
\\ &\ge   \sup_{\|g\|_{X'}=1} \int_{\Lambda} \bigg(\int_{-1}^1 |g(t)|  
 \bigg| \sum_{n=1}^N\sigma_n a_nT\left(\chi_{A_n}\right)(t)\bigg|\,dt\bigg)d\sigma
 \nonumber
\\ & =  \sup_{\|g\|_{X'}=1} \int_{-1}^1 |g(t)|  \bigg( \int_{\Lambda}
 \bigg| \sum_{n=1}^N\sigma_n a_nT\left(\chi_{A_n}\right)(t)\bigg|\,d\sigma\bigg)dt.
\end{align}

Consider now the inner integral over $\Lambda$ 
in the last term \eqref{2} of the 
previous expression. For $t\in(-1,1)$ fixed,  set 
$$
\beta_n:=a_nT\left(\chi_{A_n}\right)(t),\quad n=1.\dots,N.
$$ 
It is known  that the coordinate projections 
$$
P_n:\sigma\in\Lambda\mapsto \sigma_n\in\{-1,1\}, \quad n=1,\dots,N,
$$
form an orthonormal set, that is,
\begin{equation*}\label{orto}
\int_\Lambda P_jP_k\,d\sigma=\int_\Lambda \sigma_j\sigma_k\,d\sigma
=\delta_{j,k},\quad j,k=1,\dots,N.
\end{equation*}
Then, for the inner integral in \eqref{2}, we have
$$
 \int_{\Lambda}  \bigg| \sum_{n=1}^N\sigma_n a_nT\left(\chi_{A_n}\right)(t)\bigg|\,d\sigma
=
 \int_{\Lambda}  \bigg| \sum_{n=1}^N\beta_nP_n(\sigma)\bigg|\,d\sigma.
$$
Apply the Khintchine inequality, 
\cite[Inequality 1.10 and p.23]{diestel-jarchow-tonge}, 
 for $\{P_n\}_{n=1}^N$ yields
$$
 \int_{\Lambda}  \bigg| \sum_{n=1}^N\beta_nP_n(\sigma)\bigg|\,d\sigma
 \ge  \frac{1}{\sqrt2}   
 \bigg( \sum_{n=1}^N |\beta_n|^2\bigg)^{1/2}.
$$
Accordingly,
\begin{equation}\label{3}
 \int_{\Lambda}  \bigg| \sum_{n=1}^N\sigma_n a_nT\left(\chi_{A_n}\right)(t)\bigg|\,d\sigma
 \ge \frac{1}{\sqrt2}    
 \bigg( \sum_{n=1}^N |a_n|^2\left|T\left(\chi_{A_n}\right)(t)\right|^2\bigg)^{1/2}.
\end{equation}
Then, from \eqref{2} and \eqref{3}, it follows that
\begin{align}\label{4}
\|F\|_{L^\infty(\Lambda)}
 & \ge \frac{1}{\sqrt2}  \sup_{\|g\|_{X'}=1} \int_{-1}^1 |g(t)|   
 \bigg( \sum_{n=1}^N |a_n|^2\left|T\left(\chi_{A_n}\right)(t)\right|^2\bigg)^{1/2}\,dt\nonumber
\\ & = \frac{1}{\sqrt2}  \bigg\|   
 \bigg( \sum_{n=1}^N |a_n|^2\left|T\left(\chi_{A_n}\right)\right|^2\bigg)^{1/2}\bigg\|_X.
\end{align}

We recall the following consequence of the Stein-Weiss formula for
the distribution function of the Hilbert transform $H$ on $\R$ 
of a characteristic function, due to Laeng, 
\cite[Theorem 1.2]{laeng}. Namely, for $A\subseteq\mathbb{R}$ with 
$m(A)<\infty$ (where $m$ also denotes Lebesgue measure in $\R$), we have
$$
m(\{x\in A:\left| H(\chi_A)(x))\right|>\lambda\}) = 
\frac{2m(A)}{e^{\pi\lambda}+1},\quad \lambda>0.
$$
In particular, for any set $A\subseteq(-1,1)$ it follows, for each $\lambda>0$, that
\begin{align*}
m(\{x\in A:\left| T(\chi_A)(x)\right|>\lambda\})
= 
m(\{x\in A:\left| H(\chi_A)(x)\right|>\lambda\})
= 
\frac{2m(A)}{e^{\pi\lambda}+1}.
\end{align*}
That is,
\begin{equation}\label{extra}
m(\{x\in A:\left| T(\chi_A)(x)\right|>\lambda\})
=
\frac{2m(A)}{e^{\pi\lambda}+1},
\quad
A\in\mathcal{B},\;\; \lambda>0.
\end{equation}
Set $\lambda=1$ and $\delta:=2/(e^{\pi}+1)<1$. For each $n=1,\dots,N$, define
\begin{equation*}
A_n^1:=\{x\in A_n:\left| T(\chi_{A_n})(x)\right|>1\}.
\end{equation*}
Then  \eqref{extra} implies that
\begin{equation}\label{e}
m(A_n^1)= \frac{2m(A_n)}{e^{\pi}+1}=\delta m(A_n), \quad n=1,\dots,N.
\end{equation}
Since the sets $A_1,\dots, A_N$ are pairwise disjoint,  
so are their subsets $A_1^1,\dots, A_N^1$. Note that
$|T\left(\chi_{A_n}\right)(x)|>1$ for $x\in A_n^1$, for $n=1,\dots,N$. Thus,
on $(-1,1)$ we have the pointwise estimates
\begin{align}\label{extra2}
\bigg( \sum_{n=1}^N |a_n|^2\left|T\left(\chi_{A_n}\right)\right|^2\bigg)^{1/2}
&\ge 
\bigg( \sum_{n=1}^N |a_n|^2\left|T\left(\chi_{A_n}\right)\right|^2\chi_{A_n}\bigg)^{1/2}\nonumber
\\ & =
\sum_{n=1}^N |a_n|\left|T\left(\chi_{A_n}\right)\right|\chi_{A_n} \nonumber
\\ &\ge 
\sum_{n=1}^N |a_n|\chi_{A_n^1}.
\end{align}
Since $\|\cdot\|_X$ is a lattice norm, \eqref{extra2} yields
\begin{equation}\label{7}
\bigg\| \bigg( \sum_{n=1}^N |a_n|^2\left|T\left(\chi_{A_n}\right)\right|^2\bigg)^{1/2}\bigg\|_X
\ge 
\bigg\| \sum_{n=1}^N |a_n|\chi_{A_n^1}\bigg\|_X = \|\varphi\|_X,
\end{equation}
where $\varphi$  is the simple function
$$
\varphi:=\sum_{n=1}^N a_n\chi_{A_n^1}.
$$
From  \eqref{e}  it follows that 
\begin{equation}\label{extra4}
m(\{x\in(-1,1): |\varphi(x)|>\lambda\})
= \delta m(\{x\in(-1,1): |\phi(x)|>\lambda\}),
\quad \lambda>0.
\end{equation}
Consider the dilation operator 
$E_{\delta}\colon \widetilde X\to \widetilde X$ for   $\delta<1$ given above,
 that is, $E_{\delta}(f)(t)=f(\delta t)$ for  $0\le s\le \min\{2,1/\delta\}$ and zero otherwise.
For  the decreasing rearrangements $ \phi^*$ and $\varphi^*$ of 
 $\phi$ and $\varphi$, respectively, it follows from \eqref{extra4}  that 
$$
\phi^*=E_{\delta}(\varphi^*).
$$ 
Consequently, with $\|E_{\delta}\|$ denoting the operator norm of
$E_\delta\colon \widetilde X\to\widetilde X$, we have
\begin{align}\label{8}
\|\phi\|_X=\|\phi^*\|_{\widetilde X}&
=\|E_{\delta}(\varphi^*)\|_{\widetilde X}\le \|E_{\delta}\|
\cdot \|\varphi^*\|_{\widetilde X}
=\|E_{\delta}\|\cdot \|\varphi\|_X.
\end{align}
It follows, from \eqref{1}, \eqref{4}, \eqref{7} and \eqref{8} that
\begin{align*}
\|\phi\|_X &  \le \|E_{\delta}\| \cdot \|\varphi\|_X 
\\ & \le 
\|E_{\delta}\| \cdot
\bigg\| \bigg( \sum_{n=1}^N |a_n|^2\left|T\left(\chi_{A_n}\right)\right|^2\bigg)^{1/2}
\bigg\|_X
\\ & 
\le \sqrt2 \|E_{\delta}\| \cdot  \|\phi\|_{[T,X]}.
\end{align*}
That is, there exists a constant $M>0$, depending exclusively on $X$, 
such that
\begin{align}\label{9}
M \|\phi\|_X \le\|\phi\|_{[T,X]},
\end{align}
for \textit{all} simple functions $\phi$.

In order to extend \eqref{9} to all functions in $[T,X]$  fix $f\in[T,X]$. 
For every simple function $\phi$ satisfying $|\phi|\le |f|$  
it follows from \eqref{9} that
$$
M\|\phi\|_X \le \|\phi\|_{[T,X]} \le\|f\|_{[T,X]}.
$$
Taking the supremum with respect to all such  $\phi$ yields, 
via the Fatou property of $X$, that $f\in X$ and
$$
M\|f\|_X \le  \|f\|_{[T,X]}.
$$
In particular, $[T,X]\subseteq X$. Consequently,  $[T,X]=X$ with equivalent norms. Thus,
no genuine  $X$-valued extension of $T_X\colon X\to X$ is possible. \qed


The above Theorem has an immediate consequence. Namely, it extends, 
to \textit{all}  r.i.\ spaces with non-trivial Boyd indices, certain results known for 
those r.i.\ spaces $X$ satisfying $0<\underline{\alpha}_X\le\overline{\alpha}_X<1/2$ or
$1/2<\underline{\alpha}_X\le\overline{\alpha}_X<1$, \cite[Corollary 4.8]{curbera-okada-ricker},
and for $X=L^2(-1,1)$, \cite[Corollaries 5.4 and 5.5]{curbera-okada-ricker}.

\begin{nonumbercorollary}
Let $X$ be a r.i.\ space on $(-1,1)$ with non-trivial Boyd indices. 
\begin{itemize}
\item[(a)] There exists a constant $\beta>0$ such that, for every $f\in X$, we have 
$$
\frac{\beta}{4}\|f\|_X
\le  \sup_{A\in\mathcal{B}}\big\|T_X(\chi_A f)\big\|_X
\le \sup_{|\theta|=1}\big\|T_X(\theta f)\big\|_X
\le \sup_{|h|\le|f|}\big\|T_X(h)\big\|_X
\le \|f\|_X .
$$
\item[(b)]
For a function $f\in L^1$ the following conditions are equivalent.
\begin{itemize}
\item[(a)] $f\in X$.
\item[(b)] $T(f\chi_A)\in X$ for every $A\in\mathcal{B}$.
\item[(c)] $T(f\theta)\in X$ for every  $\theta\in L^\infty$ with $|\theta|=1$ a.e.
\item[(d)] $T(h)\in X$ for every $h\in L^0$ with $|h|\le |f|$ a.e.
\end{itemize}
\end{itemize}
\end{nonumbercorollary}


\begin{acknowledgement}
We are grateful to the referee for some valuable suggestions.
\end{acknowledgement}



\end{document}